\documentclass{article} 
\usepackage{iclr2021_conference,times}

\iclrfinalcopy

\usepackage{amsmath,amsfonts,bm}

\def\eqref#1{equation~\ref{#1}}

\def\1{\bm{1}}

\DeclareMathAlphabet{\mathsfit}{\encodingdefault}{\sfdefault}{m}{sl}
\SetMathAlphabet{\mathsfit}{bold}{\encodingdefault}{\sfdefault}{bx}{n}

\usepackage{hyperref}
\usepackage{url}
\usepackage[utf8]{inputenc} 
\usepackage[T1]{fontenc}    
\usepackage{booktabs}       
\usepackage{amsfonts}       
\usepackage{nicefrac}       
\usepackage{microtype}      
\usepackage{graphicx}
\usepackage{amsmath}
\usepackage[inline]{enumitem}
\usepackage{booktabs}
\usepackage{appendix}
\usepackage{subcaption}

\title{Data-driven Optimization Model for Global Covid-19 Intervention Plans}
\author{Chang Liu \& Akshay Budhkar \\
Georgian \\
\texttt{\{chang,akshay\}@georgian.io} \\
}

\begin{document}

\maketitle

\begin{abstract}
In the wake of COVID-19, every government huddles to find the best interventions that will reduce the number of infection cases while minimizing the economic impact. However, with many intervention policies available, how should one decide which policy is the best course of action? In this work, we describe an integer programming approach to prescribe intervention plans that optimizes for both the minimal number of daily new cases and economic impact. We present a method to estimate the impact of intervention plans on the number of cases based on historical data. Finally, we demonstrate visualizations and summaries of our empirical analyses on the performance of our model with varying parameters compared to two sets of heuristics.
\end{abstract} 

\section{Introduction}
As the pandemic continues to rage, developing data-driven approaches to predict future infection cases and prescribe actionable recommendations to prevent further infections becomes imperative. Thanks to the large-scale and recent data curation related to the pandemic, many insights can be taken from historical trends to understand the infection rate and prevention plans' effectiveness. Based on these insights, we can help decision-makers understand the impact of different intervention plans (IP) and prescribe the most beneficial policies. In order to tackle this problem, we set out to formulate an integer program, with parameters computed from historical data, to find the optimal intervention policies. The final section of this paper highlights the performance of our proposed model compared with existing heuristics. This work results from participating in the XPRIZE Pandemic Response Challenge \citep{miikkulainen:arxiv20npi}. We will be using the same dataset and standard predictor as provided in the challenge.

\vspace{-0.3cm}
\paragraph{Data}
We use the Oxford dataset \citep{hale2020variation} which contains case numbers and historical IPs taken across 280 regions around the globe since the rise of the pandemic. We can see a sample of the historical IPs prescribed in the US and Canada in Appendix \ref{IP_history}.

\paragraph{Standard Predictor}
In this study, we make use of the predictor made available by XPRIZE \citep{miikkulainen:arxiv20npi} to help us estimate parameters and make future prescriptions. The authors use a 21-day lookback window to capture the number of cases and the IPs separately. The model enforces monotonicity on the IPs as increasing their stringency should lead to fewer infection cases. They also use a 7-day rolling window to average out any inconsistencies in daily case counts. The authors then train a model comprising two separate recurrent LSTM layers as two distinct pathways that eventually lead to a single prediction. The standard predictor uses the Oxford dataset \citep{hale2020variation} until Jan 11, 2021.

\paragraph{Problem Definition}
Let $\mathcal{P}$ represent the set of IPs, in this problem, $|\mathcal{P}| = 12$, and $\mathcal{L}_p$ represent the set of available restriction level for plan $p \in \mathcal{P}$. The list of intervention plans and their associated restriction level is shown in Appendix \ref{IP}. For each region, given the different IPs and the restriction levels for each IP, the problem is to find the optimal set of IPs that minimizes both the daily new cases and the stringency cost for each region.


\section{Exact Approach to Prescribing Optimal Intervention Plans}

We present an exact algorithm to solve the task of prescribing optimal intervention plans. To formulate the problem as an integer program, we must first estimate several parameters. The formulation uses the following variables: 

\begin{itemize} 
    \item[] Let $x_{p,l} =$
    $\begin{cases}
      1 & \text{if IP $p$ at level $l$ is prescribed}\\
      0 & \text{otherwise}
    \end{cases}$    
\end{itemize}

\subsection{Stringency Costs}
Let $S_{p,l}$ represent the stringency cost associated with IP $p$ at level $l$. In this work, we experimented with three sets of stringency costs: fixed, random, and realistic. For all of these costs, we only generate a cost for each IP at level one. We multiply the cost at level one with the current level to reflect increasing costs with increasing restrictions for each increasing level. The fixed costs assume that the stringency cost for each IP is 1. We generate the random costs at random. The realistic costs are chosen based on the relative costs of implementing these IPs in Canada on a scale of 1 to 10 (10 being the highest cost), as shown in Appendix \ref{IP}. We then apply the same costs to all the regions and expect it to add value given its relatively horizontal and realistic application.

\subsection{Impact of Intervention Plans on the Number of Infection Cases}
Let $C_{p,l}$ represent the impact IP $p$ at level $l$ has on the number of infection cases. We estimate $C_{p,l}$ as follows. For each region, first, we use the standard predictor to compute baseline infection cases ($B$) by setting all IPs to level zero, simulating a scenario where a region prescribes no IP. Next, for each IP and each restriction level, we use the standard predictor to compute the estimated number of infection cases ($E_{p,l}$) by only activating that specific IP at that restriction level. For each region, we then compute $C_{p,l}$ for each IP and its associated level as shown in Equation \ref{weights_eq}.

\begin{equation} \label{weights_eq}
    C_{p,l} = \frac{(E_{p,l} - B)}{B} * 100
\end{equation}

To test the impact of these weights, we perform this estimating process for one and seven days in the future, starting from three different days (Aug 2, 2020; Jan 2, 2021; Jan 15, 2021), creating in total six sets of impacts that we test with our optimization model. 

\subsection{Objective Function Normalization}

Since we are optimizing for two different objectives: the number of cases, and stringency costs, we need to normalize the two numbers so the optimization model does not lean towards one or the other. We first normalize the stringency costs at level one for the 12 IPs to sum to 1. The total stringency cost can be written as in Objective \ref{stringency}. Since the restriction level for any IP is at most 4, the total stringency cost is constrained to be at most 4. 

\begin{equation} \label{stringency}
    \displaystyle \sum_{\substack{p \in \mathcal{P}}} \displaystyle \sum_{\substack{l \in \mathcal{L}_p}} S_{p, l} * x_{p, l}
\end{equation}

As for the number of cases, the impacts are at most 0.1 and on average at 0.02. We then normalize the number of infections to the current case numbers. For each region, let \(\beta\) be the number of current new cases. Let \(\alpha\) be the number of new cases in the previous day, which we compute using the standard predictor. The number of new cases can be written as in Objection \ref{cases}. As we prescribe for multiple days, \(\beta\) remains the same while \(\alpha\) will vary. As \(\alpha\), the number of predicted future new cases increases, Objective \ref{cases} will increase, forcing the model to pay more weights to the number of cases. On the contrary, if \(\alpha\) decreases, the model will start to shift more weights to the stringency costs.

\begin{equation} \label{cases}
    \frac{1}{\beta} \left( \alpha + \displaystyle \sum_{\substack{p \in \mathcal{P}}} \displaystyle \sum_{\substack{l \in \mathcal{L}_p}} C_{p, l} * x_{p, l} * \alpha \right)
\end{equation}

\subsection{Integer Program Formulation}
To prescribe for several days and multiple regions, we take an iterative approach. For each prescription day \(d\) and for each region \(r\), we solve the integer program shown in Figure \ref{integer_program}. This formulation uses binary variables \(x_{p,l} = 1\) if IP \(p\) at level \(l\) is prescribed. The objective function \ref{obj} minimizes the total number of new cases and stringency costs. In Equation \ref{c2}, we make sure that for each IP, we select only one restriction level at a time. 

\vspace{-0.3cm}
\begin{figure}[htbp]
    \centering
\begin{align}
    \mbox{minimize   } &  \text{Objectives } (\ref{stringency}) + (\ref{cases}) \label{obj}\\
    \mbox{s.t.   } 
    & \displaystyle \sum_{\substack{l \in \mathcal{L}_p}} x_{p,l} = 1  & \forall p \in \mathcal{P} \label{c2}\\
    & x_{p,l} \in \{0,1\} & \forall p \in \mathcal{P}, l \in \mathcal{L}_p \label{c1}
\end{align}
\caption{Integer Program Model}
\label{integer_program}
\end{figure}
\vspace{-0.3cm}

\subsection{Consecutive Constraint}
In the real world, it is natural that when an IP is put in place, it will last for at least several days. By analyzing historical IPs across all the regions (Appendix \ref{IP_history}), we have created a minimum number of consecutive days that each intervention plan at a specific level must hold, which we set a ceiling of 7 days, denoted as \(D_{p,l}\). At the time of prescribing for a day, based on the previous seven days of prescribed IPs, we create a new parameter \(F_{p,l} = 1\) if IP \(p\) at level \(l\) has not been prescribed to enough consecutive days, thus forcing the model to set \(x_{p,l} = 1\), and 0 otherwise. We formulate this new constraint as Equation \ref{consecutive} and add it to the integer program formulation.

\vspace{-0.3cm}
\begin{align}
    & D_{p,l} \leq x_{p,l} & \forall p \in \mathcal{P}, l \in \mathcal{L}_p \label{consecutive}
\end{align}

\section{Experimental Results}
In this section, we present our integer program's performance, with and without the consecutive constraint, and with varying weights, compared with two sets of heuristics provided by XPRIZE\footnote{\url{https://github.com/leaf-ai/covid-xprize/}}. 

\vspace{-0.3cm}
\paragraph{Blind-greedy Heuristic}
The blind-greedy starts by setting all IPs to level zero, then iteratively, increases a restriction level of the next least costly IP. In this way, each subsequent prescription sets higher levels of restrictions but at the same time also becomes more costly. 

\vspace{-0.3cm}
\paragraph{Random Heuristic}
For each region and each prescription day, this heuristic assigns a valid level for each IP at random. 

We make prescriptions for 28 days from Jan 12, 2021. We then measure the mean stringency and the mean number of new cases predicted by the standard predictor. The optimization solver used is the SCIP solver \citep{GamrathEtal2020OO} imported from Google OR Tools \citep{ortools}. We set the solver to use its default parameters.

\begin{figure}[htbp]
    \centering
    \begin{subfigure}[t]{0.4\textwidth}
        \centering
        \includegraphics[height=1.7in]{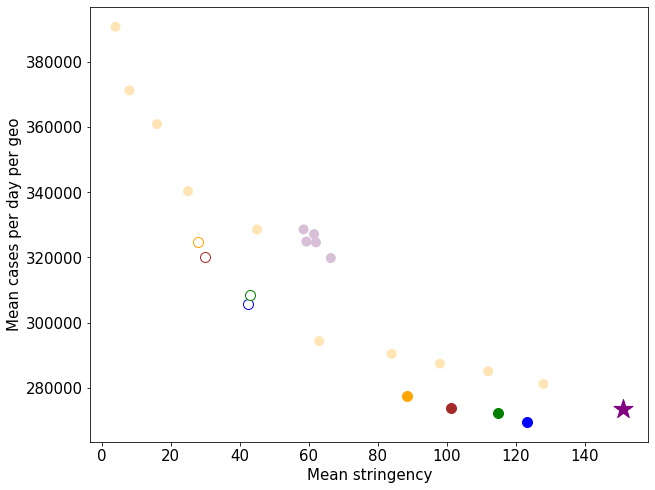}
        \caption{United States}
    \end{subfigure}%
    ~ 
    \begin{subfigure}[t]{0.6\textwidth}
        \centering
        \includegraphics[height=1.7in]{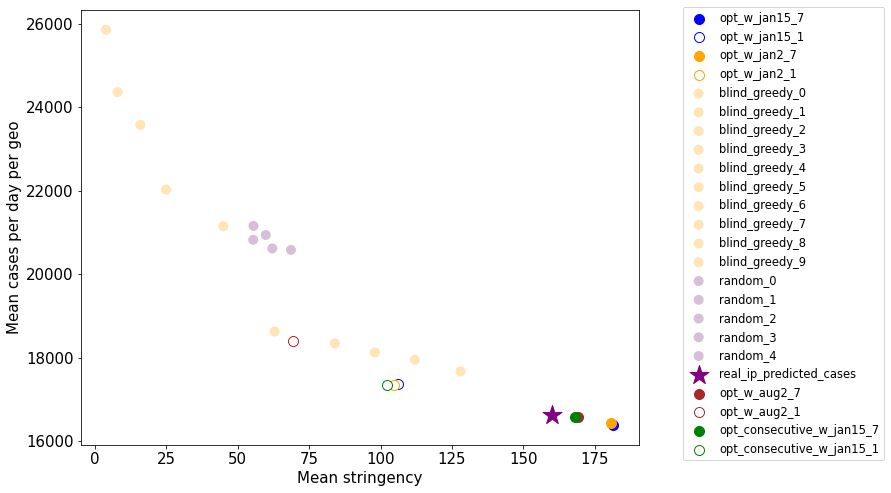}
        \caption{Canada}
    \end{subfigure}
    \caption{Performance of different models with realistic stringency costs in the US and Canada}
    \label{result}
\end{figure}

In Figure \ref{result}, we can see the performance of various models using the realistic stringency costs in the US and Canada. As expected, the random models perform the worst. The blind\_greedy models generate a Pareto front highlighting the trade-off between the two factors at hand. We can observe that the real IPs tend to be on the more expensive side for both the US and Canada. A few of our optimization models' weights find less costly solutions for the US while leading to fewer cases. In Canada, some of our solutions trade-off a few extra cases to reduce the stringency impact significantly. Most of the real IPs and our optimization solutions do better than our baseline heuristics. We report the detailed numbers in Appendix \ref{results_all}.

In Appendix \ref{Predicted_IP}, we compare the actual IPs that we incorporated by the US and Canada with a couple of our models. While our models' performance with and without the consecutive constraint is comparable, the nonconsecutive prescriptions are unrealistic given their oscillating nature, and it is not as helpful given the probable cost of turning these IPs on and off in real life. Our seven-day, Jan 15 weights lead to prescriptions similar to the US's actual interventions, but the slightly stricter public transportation rules could have led to a better outcome in both the stringency and the number of cases. For Canada, the one-day, Jan 15 weights lead to a slight increase in the number of cases, but laxer stay-at-home requirements coupled with stricter testing policies could significantly reduce the economic impact.


Furthermore, we have designed live interactive dashboards that allow users to study all the different models' efficacy across regions and different sets of stringency costs. The \hyperlink{https://tabsoft.co/3ejM0iq}{first dashboard}\footnote{\url{https://tabsoft.co/3ejM0iq}} compares the trade-off between the number of cases and stringency costs as seen in Figure \ref{fig1} in Appendix \ref{dashboard}. The \hyperlink{https://tabsoft.co/3d71wOd}{second dashboard}\footnote{\url{https://tabsoft.co/3d71wOd}} simulates the case numbers based on the different IPs over a 28-day period as seen in Figure \ref{fig2} in Appendix \ref{dashboard}. Users can choose different regions, weight types, models and study these factors' effect on the case numbers over 28 days starting Jan 12, 2021.


These dashboards have an immediate short-term application to the stakeholders dealing with public health decisions as economies start to re-open. For example, based on the first dashboard with Jan 15 weights prescriptions, the Canadian government could reduce the average economic impact by approximately 22\% for a slight increase in the number of cases over a 28-day horizon.

\section{Conclusions and Future Work}
In this work, we have described a data-driven integer program approach to optimal intervention plan prescription. We presented the performance of our models in comparison with other heuristics and introduce a framework to visualize the impact of intervention plans on the number of cases and economic costs on a region-by-region basis.

Relevant stakeholders in public health, epidemiological research, and governments can readily implement the recommendations to manage the spread of the pandemic while minimizing the impact on society. We look forward to presenting our framework to our local government to get feedback and to get to a broader audience through this platform to save lives while mitigating societal impact. 

This work is preliminary and there are multiple ways we look forward to extending our work. First, we believe we can improve on the normalization of the two costs in the integer program objective function by implementing algorithms from the multi-objective optimization field. Further, we can also improve on the estimated \(C_{p,l}\) by performing a linear regression on the historical IPs. We have also started exploring Reinforcement Learning methods to tackle this problem, where we first model the spread of the disease as an environment for an agent to explore, we then integrate Deep Q-Learning \citep{DQN} and Soft Actor Critic (SAC) algorithms \citep{haarnoja2018soft} to optimize IPs for each region. We believe this might be a promising approach and will investigate further. Finally, the performance of the prescription model is heavily dependent on the performance of the infection cases predictor. We believe in our future work, we can also improve on the standard predictor or incorporate the two models in an interactive system \citep{elmachtoub2017smart}.

\newpage

\bibliographystyle{iclr2021_conference}
\bibliography{ref}

\newpage

\appendix

\section{Historical Intervention Plans in the US and Canada } \label{IP_history}

\begin{figure}[htbp]
    \centering
    \includegraphics[width=0.9\textwidth]{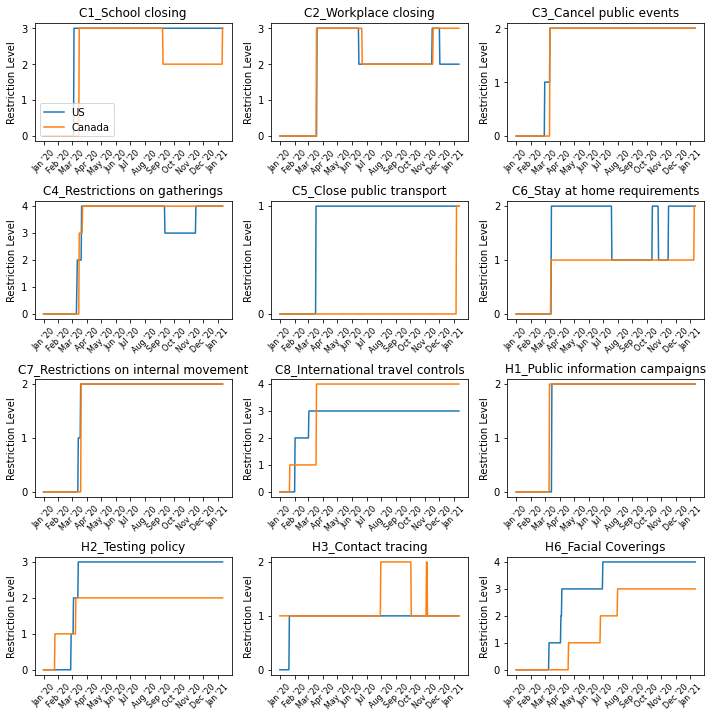}
    \caption{Trends of IPs in the US and Canada}
    \label{IP_history_figure}
\end{figure}

\section{Intervention Plans} \label{IP}

\begin{table}[htbp]
\centering
\begin{tabular}{lrr}
\toprule
\textbf{Intervention Plan} & \textbf{Restriction Levels} & \textbf{Realistic Stringency Cost} \\
\midrule
C1\_School closing & (0,1,2,3) &	9 \\
C2\_Workplace closing &(0,1,2,3) &	6 \\
C3\_Cancel public events &(0,1,2) & 	2 \\
C4\_Restrictions on gatherings & (0,1,2,3,4) &	5 \\
C5\_Close public transport & (0,1,2) &	8 \\
C6\_Stay at home requirements & (0,1,2,3) &	7 \\
C7\_Restrictions on internal movement & (0,1,2) &	7 \\
C8\_International travel controls &	(0,1,2,3,4)&8 \\
H1\_Public information campaigns &(0,1,2)&	2 \\
H2\_Testing policy &(0,1,2,3)&	3 \\
H3\_Contact tracing &(0,1,2)&	7 \\
H6\_Facial Coverings &(0,1,2,3,4)&	2 \\
\bottomrule
\end{tabular}
\caption{Restriction levels and realistic stringency cost of IPs}
\label{table1}
\end{table}

\section{Dashboards} \label{dashboard}
\begin{figure}[htbp] 
    \centering
    \includegraphics[width=0.70\textwidth]{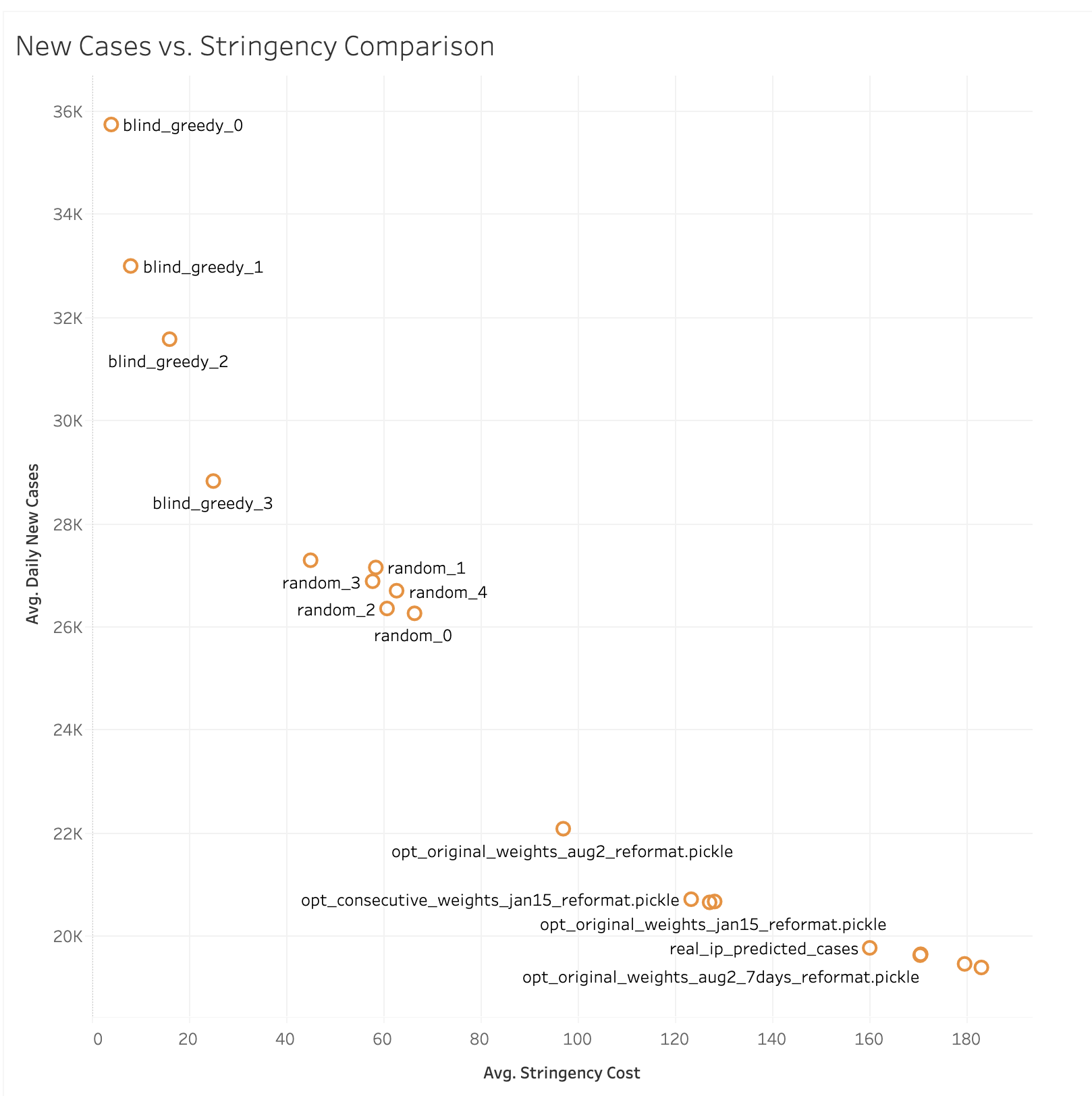}
    \caption{Comparison of daily average cases and Stringency of different models}
    \label{fig1}
\end{figure}

\begin{figure}[htbp] 
    \centering
    \includegraphics[width=0.75\textwidth]{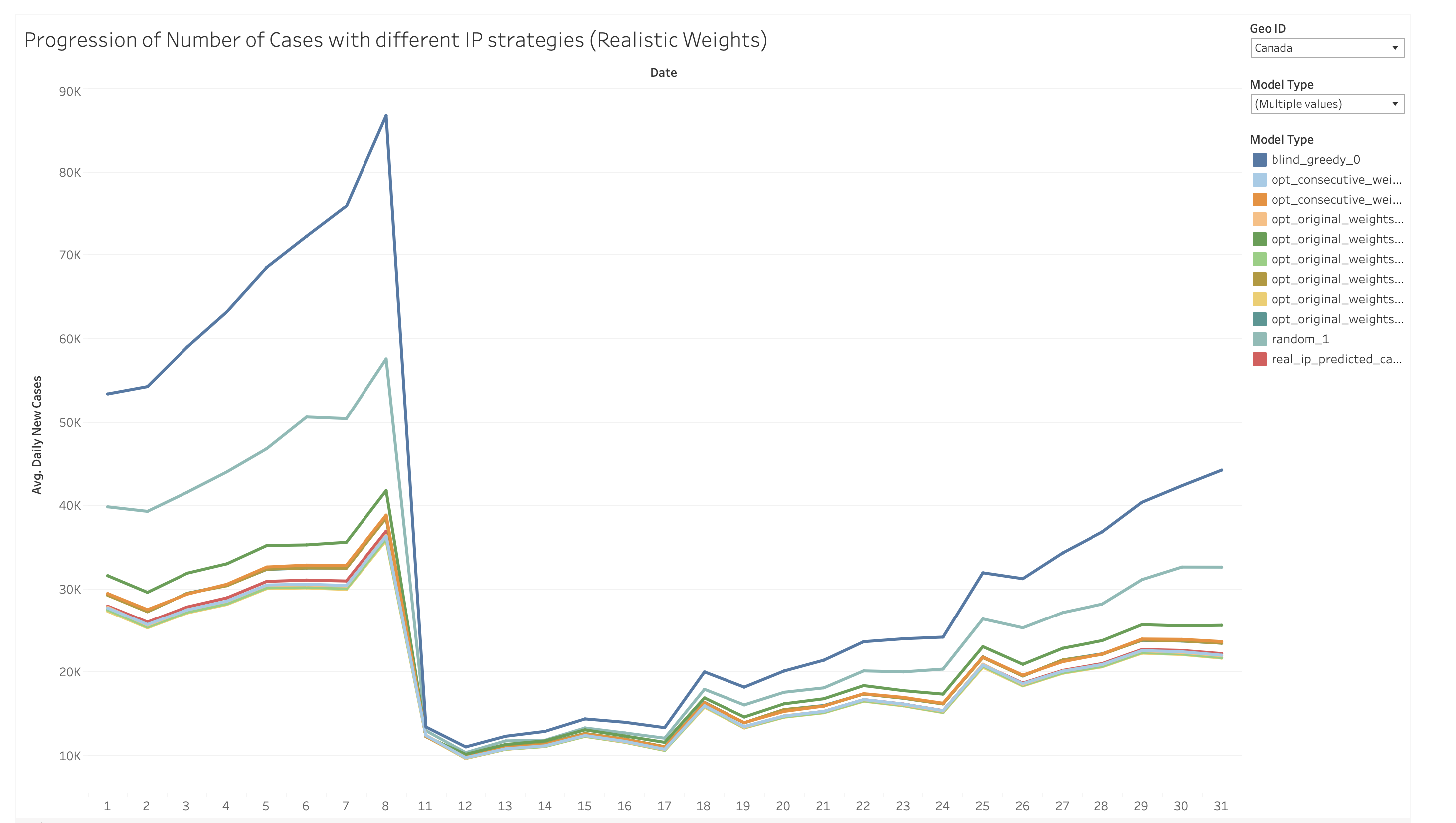}
    \caption{Impact of different IPs over the a seven day horizon}
    \label{fig2}
\end{figure}

The \hyperlink{https://tabsoft.co/3ejM0iq}{first dashboard}\footnote{\url{https://tabsoft.co/3ejM0iq}} compares the trade off between number of cases and stringency costs, a snapshot is shown in Figure \ref{fig1}. The \hyperlink{https://tabsoft.co/3d71wOd}{second dashboard} \footnote{\url{https://tabsoft.co/3d71wOd}} simulates the case numbers based on the different IPs over a 28-day period, a snapshot is shown in Figure \ref{fig2}. 


\newpage

\section{Performance Results} \label{results_all}
We present the detailed number of our experimental results here. There are in total 10 blind-greedy models with 5 random models. We experimented with 6 sets of different \(C_{p,l}\), predicting 1 day and 7 days starting from 3 different dates: Aug 2, 2020, Jan 2, 2021; and Jan 15, 2021. We further experimented the consecutive constraint with two sets of \(C_{p,l}\). 
\subsection{United States}
\renewcommand{\arraystretch}{1.5}
\begin{table}[htbp]
\begin{tabular}{|l|r|r|r|r|r|r|}
\hline
\textbf{Stringency Costs}& \multicolumn{2}{c|}{\textbf{fixed}} & \multicolumn{2}{c|}{\textbf{random}} & \multicolumn{2}{c|}{\textbf{realistic}} \\
\hline
\textbf{Measure}& \multicolumn{1}{c|}{\textbf{cases}}	&\multicolumn{1}{c|}{\textbf{costs}}	&\multicolumn{1}{c|}{\textbf{cases}}	&\multicolumn{1}{c|}{\textbf{costs}}	&\multicolumn{1}{c|}{\textbf{cases}}	&\multicolumn{1}{c|}{\textbf{costs}} \\
\hline
opt\_consecutive\_w\_jan15\_7 & 273939.20          & 21.33          & 269403.72          & 19.95          & 272183.34          & 114.81          \\ \hline
opt\_consecutive\_w\_jan15\_1 & 307375.74          & 7.71           & 307685.89          & 8.58           & 308403.71          & 42.95           \\ \hline
opt\_w\_jan15\_7              & 268771.90          & 24.24          & 268738.99          & 19.97          & 269436.13          & 123.10          \\ \hline
opt\_w\_jan15\_1              & 300796.03          & 8.67           & 300927.83          & 8.09           & 305727.94          & 42.43           \\ \hline
opt\_w\_jan2\_7               & 276585.64          & 17.43          & 272252.66          & 17.09          & 277367.53          & 88.38           \\ \hline
opt\_w\_jan2\_1               & 326332.17          & 5.10           & 322437.11          & 4.03           & 324601.12          & 28.00           \\ \hline
opt\_w\_aug2\_7               & 273110.33          & 20.48          & 270715.27          & 18.18          & 273755.98          & 101.14          \\ \hline
opt\_w\_aug2\_1               & 318821.35          & 5.67           & 317196.32          & 4.70           & 320039.63          & 29.90           \\ \hline
blind\_greedy\_0              & 350768.63          & 3.00           & 276585.64          & 17.43          & 390854.80          & 4.00            \\ \hline
blind\_greedy\_1              & 306824.18          & 6.00           & 272252.66          & 17.09          & 371281.95          & 8.00            \\ \hline
blind\_greedy\_2              & 302727.20          & 8.00           & 276585.64          & 17.43          & 360932.75          & 16.00           \\ \hline
blind\_greedy\_3              & 292355.06          & 12.00          & 272252.66          & 17.09          & 340350.61          & 25.00           \\ \hline
blind\_greedy\_4              & 285647.76          & 14.00          & 293142.51          & 9.98           & 328540.60          & 45.00           \\ \hline
blind\_greedy\_5 & 281362.98 & 17.00 & 289671.00 & 11.56 & 294287.71 & 63.00  \\ \hline
blind\_greedy\_6 & 278267.20 & 19.00 & 287302.63 & 14.04 & 290403.70 & 84.00  \\ \hline
blind\_greedy\_7 & 272759.52 & 23.00 & 284184.68 & 17.82 & 287388.78 & 98.00  \\ \hline
blind\_greedy\_8 & 270096.78 & 25.00 & 280900.00 & 23.18 & 285031.79 & 112.00 \\ \hline
blind\_greedy\_9 & 265914.71 & 28.00 & 278425.29 & 25.90 & 281126.78 & 128.00 \\ \hline
random\_0                     & 322917.85          & 10.95          & 324300.98          & 9.79           & 327159.45          & 61.52           \\ \hline
random\_1                     & 328073.25          & 10.62          & 322851.06          & 11.23          & 328588.25          & 58.48           \\ \hline
random\_2                     & 324572.39          & 10.57          & 332052.66          & 10.76          & 324622.28          & 62.05           \\ \hline
random\_3                     & 330125.04          & 10.00          & 326852.32          & 10.02          & 319797.71          & 66.38           \\ \hline
random\_4                     & 326739.40          & 10.38          & 327730.99          & 10.34          & 324854.00          & 59.29           \\ \hline
real\_ip\_predicted\_cases    & 273586.82          & 29.00          & 273586.82          & 28.92          & 273586.82          & 151.00          \\ \hline
\textit{actual\_cases}        & \textit{174680.05} & \textit{29.00} & \textit{174680.05} & \textit{28.92} & \textit{174680.05} & \textit{151.00} \\ \hline

\end{tabular}
\caption{Performance of various model for prescribing 28 days from Jan 12, 2021 in the United States}
\end{table}

\newpage
\subsection{Canada}

\begin{table}[htbp]
\begin{tabular}{|l|r|r|r|r|r|r|}
\hline
\textbf{Stringency Costs}& \multicolumn{2}{c|}{\textbf{fixed}} & \multicolumn{2}{c|}{\textbf{random}} & \multicolumn{2}{c|}{\textbf{realistic}} \\
\hline
\textbf{Measure}& \multicolumn{1}{c|}{\textbf{cases}}	&\multicolumn{1}{c|}{\textbf{costs}}	&\multicolumn{1}{c|}{\textbf{cases}}	&\multicolumn{1}{c|}{\textbf{costs}}	&\multicolumn{1}{c|}{\textbf{cases}}	&\multicolumn{1}{c|}{\textbf{costs}} \\
\hline
opt\_consecutive\_w\_jan15\_7 & 16525.46          & 31.19          & 16789.03          & 27.48          & 16588.19          & 168.00         \\ \hline
opt\_consecutive\_w\_jan15\_1 & 17273.91          & 19.57          & 17450.57          & 14.24          & 17346.64          & 102.38         \\ \hline
opt\_w\_jan15\_7              & 16410.78          & 32.57          & 16480.23          & 31.43          & 16383.82          & 181.43         \\ \hline
opt\_w\_jan15\_1              & 17239.80          & 20.43          & 17451.10          & 14.11          & 17356.17          & 106.00         \\ \hline
opt\_w\_jan2\_7               & 16476.53          & 32.24          & 16556.99          & 31.07          & 16438.41          & 180.76         \\ \hline
opt\_w\_jan2\_1               & 17307.67          & 19.67          & 17476.61          & 13.88          & 17340.07          & 104.57         \\ \hline
opt\_w\_aug2\_7                & 16571.76          & 30.81          & 16780.44          & 27.34          & 16572.66          & 169.00         \\ \hline
opt\_w\_aug2\_1               & 18190.33          & 14.05          & 17892.12          & 8.34           & 18388.79          & 69.24          \\ \hline
blind\_greedy\_0              & 22811.46          & 3.00           & 22151.84          & 0.09           & 25858.13          & 4.00           \\ \hline
blind\_greedy\_1              & 19543.13          & 6.00           & 20593.91          & 1.47           & 24362.63          & 8.00           \\ \hline
blind\_greedy\_2              & 19242.32          & 8.00           & 20275.72          & 2.83           & 23577.58          & 16.00          \\ \hline
blind\_greedy\_3              & 18484.12          & 12.00          & 19852.64          & 4.35           & 22028.85          & 25.00          \\ \hline
blind\_greedy\_4              & 17996.32          & 14.00          & 19300.23          & 6.05           & 21147.75          & 45.00          \\ \hline
blind\_greedy\_5              & 17685.69          & 17.00          & 17710.85          & 8.81           & 18624.70          & 63.00          \\ \hline
blind\_greedy\_6              & 17461.80          & 19.00          & 17447.73          & 13.65          & 18341.68          & 84.00          \\ \hline
blind\_greedy\_7              & 17064.48          & 23.00          & 17228.14          & 16.15          & 18122.44          & 98.00          \\ \hline
blind\_greedy\_8              & 16872.91          & 25.00          & 16905.64          & 21.79          & 17951.33          & 112.00         \\ \hline
blind\_greedy\_9              & 16572.61          & 28.00          & 16745.72          & 26.02          & 17668.37          & 128.00         \\ \hline
random\_0                     & 21057.65          & 9.95           & 20743.72          & 11.14          & 20582.15          & 68.67          \\ \hline
random\_1                     & 20808.50 & 11.62 & 20705.43 & 11.53 & 21155.11 & 55.57 \\ \hline
random\_2                     & 20777.58          & 10.95       & 20742.32          & 10.97       & 20614.5           & 62.14       \\ \hline
random\_3                     & 20963.39          & 10.90       & 20952.72          & 11.15       & 20822.47          & 55.48       \\ \hline
random\_4                     & 20673.91          & 11.76       & 20774.19          & 11.00             & 20937.03          & 59.86       \\ \hline
real\_ip\_predicted\_cases    & 16615.24          & 29.00             & 16615.24          & 29.52          & 16615.24          & 160.00            \\ \hline
\textit{actual\_cases}        & \textit{5488.857} & \textit{29.00}    & \textit{5488.86} & \textit{29.52} & \textit{5488.857} & \textit{160.00}   \\ \hline

\end{tabular}
\caption{Performance of various model for prescribing 28 days from Jan 12, 2021 in Canada}
\end{table}

\newpage
\subsection{World}

\begin{table}[htbp]
\begin{tabular}{|l|r|r|r|r|r|r|}
\hline
\textbf{Stringency Costs}& \multicolumn{2}{c|}{\textbf{fixed}} & \multicolumn{2}{c|}{\textbf{random}} & \multicolumn{2}{c|}{\textbf{realistic}} \\
\hline
\textbf{Measure}& \multicolumn{1}{c|}{\textbf{cases}}	&\multicolumn{1}{c|}{\textbf{costs}}	&\multicolumn{1}{c|}{\textbf{cases}}	&\multicolumn{1}{c|}{\textbf{costs}}	&\multicolumn{1}{c|}{\textbf{cases}}	&\multicolumn{1}{c|}{\textbf{costs}} \\

\hline
opt\_consecutive\_w\_jan15\_7 & 6769.99          & 26.03          & 6769.69          & 23.67          & 6777.57          & 142.09          \\ \hline
opt\_consecutive\_w\_jan15\_1 & 7405.65          & 19.57          & 7338.26          & 17.31          & 7464.95          & 107.01          \\ \hline
opt\_w\_jan15\_7              & 6665.29          & 28.15          & 6684.85          & 25.74          & 6666.85          & 154.33          \\ \hline
opt\_w\_jan15\_1              & 7263.72          & 20.56          & 7254.68          & 18.00          & 7329.61          & 110.60          \\ \hline
opt\_w\_jan2\_7               & 6810.59          & 25.85          & 6799.07          & 23.50          & 6820.37          & 141.28          \\ \hline
opt\_w\_jan2\_1               & 7781.06          & 17.31          & 7690.06          & 15.08          & 7786.65          & 93.88           \\ \hline
opt\_w\_aug2\_7               & 6794.48          & 24.79          & 6794.49          & 22.21          & 6803.79          & 134.86          \\ \hline
opt\_w\_aug2\_1               & 7829.05          & 14.92          & 7734.96          & 13.01          & 7863.44          & 80.52           \\ \hline
blind\_greedy\_0              & 9863.78          & 3.00           & 10795.07         & 0.44           & 11533.58         & 4.00            \\ \hline
blind\_greedy\_1              & 8135.79          & 6.00           & 10029.18         & 1.34           & 10706.43         & 8.00            \\ \hline
blind\_greedy\_2              & 7980.38          & 8.00           & 9444.77          & 2.63           & 10277.08         & 16.00           \\ \hline
blind\_greedy\_3              & 7591.98          & 12.00          & 8775.29          & 4.27           & 9441.59          & 25.00           \\ \hline
blind\_greedy\_4              & 7344.57          & 14.00          & 8349.93          & 6.39           & 8973.26          & 45.00           \\ \hline
blind\_greedy\_5              & 7188.00          & 17.00          & 8017.03          & 8.89           & 7663.01          & 63.00           \\ \hline
blind\_greedy\_6              & 7075.65          & 19.00          & 7572.93          & 11.80          & 7518.95          & 84.00           \\ \hline
blind\_greedy\_7              & 6877.36          & 23.00          & 7326.67          & 15.23          & 7407.82          & 98.00           \\ \hline
blind\_greedy\_8              & 6782.21          & 25.00          & 7085.89          & 19.14          & 7321.33          & 112.00          \\ \hline
blind\_greedy\_9              & 6633.68          & 28.00          & 6901.52          & 23.47          & 7178.85          & 128.00          \\ \hline
random\_0                     & 8841.54          & 10.99          & 8797.98          & 10.95          & 8824.25          & 60.56           \\ \hline
random\_1                     & 8808.71          & 11.02          & 8764.34          & 10.99          & 8849.60          & 60.10           \\ \hline
random\_2                     & 8779.83          & 11.02          & 8829.24          & 11.07          & 8810.69          & 60.65           \\ \hline
random\_3                     & 8835.40          & 10.95          & 8786.90          & 10.99          & 8789.41          & 60.44           \\ \hline
random\_4                     & 8800.00          & 11.03          & 8792.29          & 11.03          & 8792.33          & 60.65           \\ \hline
real\_ip\_predicted\_cases    & 5536.39          & 21.64          & 5536.39          & 21.55          & 5536.39          & 111.10          \\ \hline
\textit{actual\_cases}        & \textit{3144.09} & \textit{21.64} & \textit{3144.09} & \textit{21.55} & \textit{3144.09} & \textit{111.10} \\ \hline

\end{tabular}
\caption{Performance of various model for prescribing 28 days from Jan 12, 2021 in the world}
\end{table}

\newpage
\section{Predicted and Real Intervention Plans in the US and Canada } \label{Predicted_IP}

\begin{figure}[htbp]
    \centering
    \includegraphics[width=\textwidth]{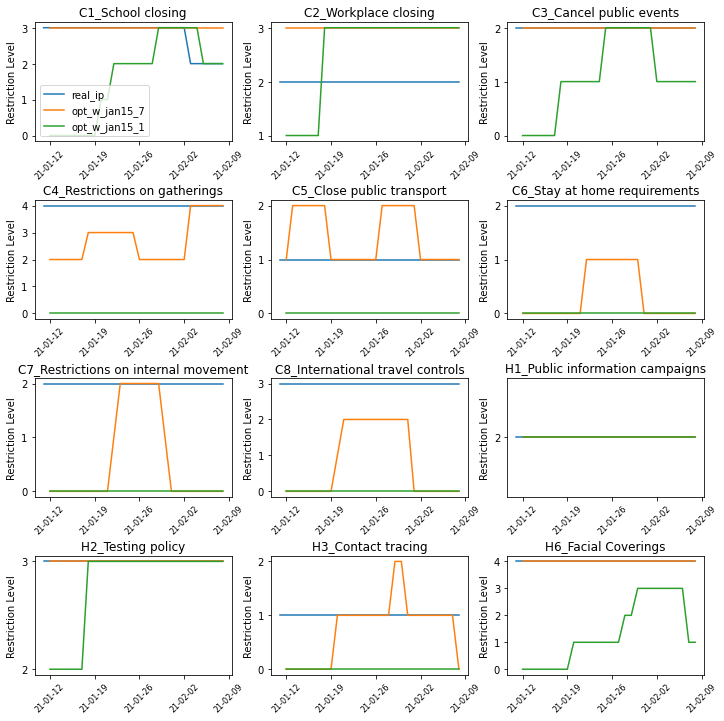}
    \caption{Predicted and real IPs with consecutive constraint in the US}
    \label{US_consecutive}
\end{figure}

\begin{figure}[htbp]
    \centering
    \includegraphics[width=\textwidth]{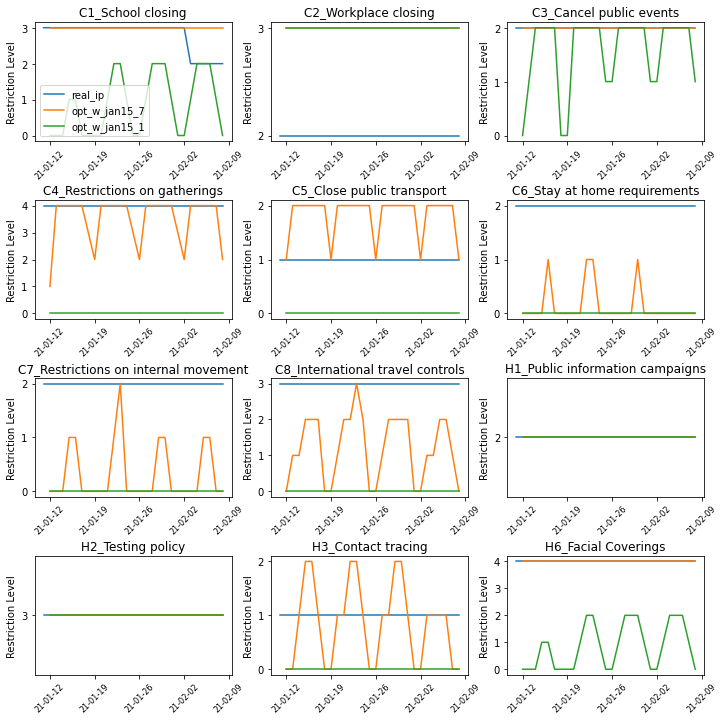}
    \caption{Predicted and real IPs in the US}
    \label{US_non_consecutive}
\end{figure}

\begin{figure}[htbp]
    \centering
    \includegraphics[width=\textwidth]{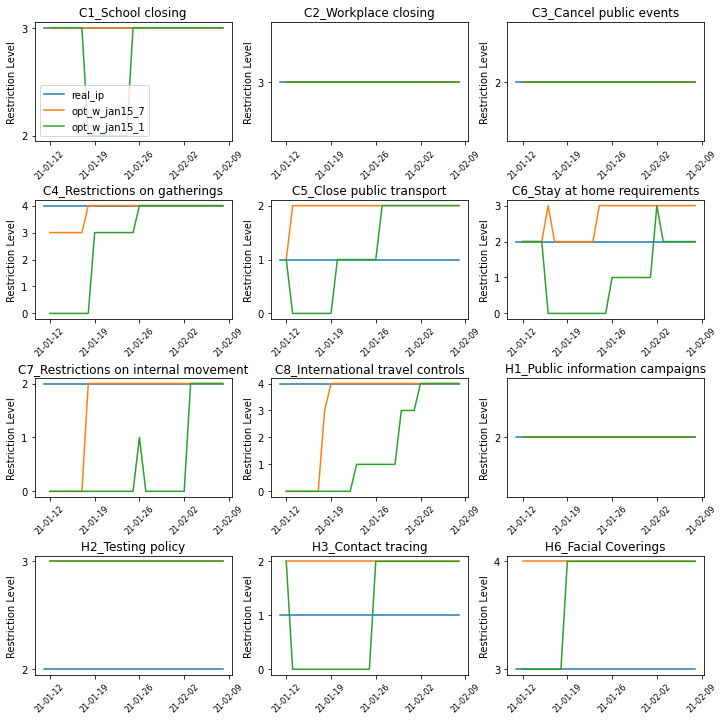}
    \caption{Predicted and real IPs with consecutive constraint in Canada}
    \label{Canada_consecutive}
\end{figure}

\begin{figure}[htbp]
    \centering
    \includegraphics[width=\textwidth]{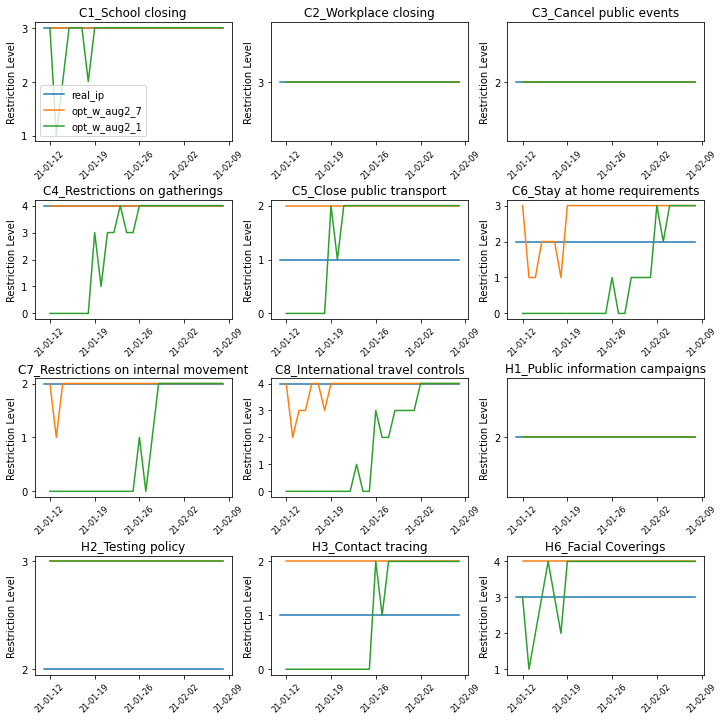}
    \caption{Predicted and real IPs in Canada}
    \label{Canada_non_consecutive}
\end{figure}

\end{document}